\documentclass[utf8]{FrontiersinHarvard} 

\usepackage{url,hyperref,lineno,microtype,subcaption}
\usepackage[onehalfspacing]{setspace}

\usepackage{tabularx}

\def\keyFont{\fontsize{8}{11}\helveticabold }
\def\firstAuthorLast{Wang {et~al.}} 
\def\Authors{Xiaoyu Wang\,$^{1,*}$ and Martin Benning\,$^{2, 3}$}
\def\Address{$^{1}$Department of Applied Mathematics and Theoretical Physics, University of Cambridge, Cambridge, United Kingdom \\
$^{2}$School of Mathematical Sciences, Faculty of Science and Engineering, Queen Mary University of London, London, United Kingdom \\
$^{3}$The Alan Turing Institute, London, United Kingdom}

\def\corrAuthor{Martin Benning}

\def\corrEmail{m.benning@qmul.ac.uk}

\newtheorem{definition}{Definition}
\newtheorem{example}{Example}
\newtheorem{lemma}{Lemma}
\newtheorem{remark}{Remark}
\newtheorem{theorem}{Theorem}

\newcommand{\actfct}{\sigma}
\DeclareMathOperator*{\argmin}{\arg \min}
\DeclareMathOperator*{\dom}{\text{dom}}
\newcommand{\bregmandis}[1][\regfct]{D_{#1}}
\newcommand{\symbreg}[1][\regfct]{\bregmandis[#1]^{\text{symm}}}
\newcommand{\bregmanloss}[1][\Psi]{B_{#1}}
\newcommand{\data}{y}

\newcommand{\jensen}[1][\Psi]{J_{#1}}
\newcommand{\layer}{L}
\newcommand{\noisebound}{\delta}
\newcommand{\noisydata}{\data^\noisebound}

\newcommand{\prox}[1][\Psi]{\text{prox}_{#1}}
\newcommand{\R}{\mathbb{R}}
\newcommand{\reconstruction}[1][\alpha]{x_{#1}}
\newcommand{\regfct}{R}
\newcommand{\regop}[1][\alpha]{\mathcal{R}_{#1}}
\newcommand{\scelem}{v^\dagger}
\newcommand{\solution}{x^\dagger}

\begin{document}
\onecolumn
\firstpage{1}

\title[Lifted Bregman inversion]{A Lifted Bregman Formulation for the Inversion of Deep Neural Networks}

\author[\firstAuthorLast ]{\Authors} 
\address{} 
\correspondance{} 

\extraAuth{}

\maketitle

\begin{abstract}

We propose a novel framework for the regularised inversion of deep neural networks. The framework is based on the authors' recent work on training feed-forward neural networks without the differentiation of activation functions. The framework lifts the parameter space into a higher dimensional space by introducing auxiliary variables, and penalises these variables with tailored Bregman distances. We propose a family of variational regularisations based on these Bregman distances, present theoretical results and support their practical application with numerical examples. In particular, we present the first convergence result (to the best of our knowledge) for the regularised inversion of a single-layer perceptron that only assumes that the solution of the inverse problem is in the range of the regularisation operator, and that shows that the regularised inverse provably converges to the true inverse if measurement errors converge to zero.

\tiny
 \keyFont{ \section{Keywords:} Inverse problems, regularisation theory, lifted network training, Bregman distance, perceptron, multi-layer perceptron, variational regularisation, total variation regularisation} 
\end{abstract}

\section{Introduction}\label{sec:introduction}



Neural networks are computing systems that have revolutionised a wide range of research domains over the past decade and outperformed many traditional machine learning approaches (cf. \citep{lecun2015deep,goodfellow2016deep}). This performance often comes at the cost of interpretability (or rather a lack thereof) of the outputs that a neural network produces for given inputs. As a consequence, a lot of research has focused on understanding representations of neural networks and on developing strategies to interpret these representations, predominantly with saliency maps \citep{simonyan2013deep,fong2017interpretable,chang2018explaining,fong2019understanding}. An alternative approach focuses on understanding deep image representations by inverting them \citep{mahendran2015understanding}. The authors propose a total-variation-based variational optimisation method that aims to infer the network input from the network output with regularised inversion. 

While the concept of inverting neural networks is certainly not new (cf. \citep{linden1989inversion,kindermann1990inversion,jensen1999inversion,lu1999inverting}), there has been increasing interest in recent years largely due to developments in nonlinear dimensionality reduction and generative modelling that include (but are not limited to) (variational) Autoencoders \citep{kingma2013auto}, Normalising Flows \citep{rezende2015variational,dinh2015nice}, Cycle-Consistent Generative Adversarial Networks \citep{zhu2017unpaired}, and Probabilistic Diffusion Models \citep{sohl2015deep,ho2020denoising}.

While several approaches for the inversion of neural networks have been proposed especially in the context of generative modelling (see for example \citep{behrmann2019invertible,behrmann2021understanding} in the context of normalising flows, \citep{xia2022gan} in the context of generative adversarial networks and \citep{gal2022image} in the context of probabilistic diffusion models), an important aspect, which is often overlooked, is that invertible operations alone are not automatically stable with respect to small variations in the data. For example, computing the solution of the heat equation after a fixed termination time is stable with respect to variations in the initial condition, but estimating the initial condition from the terminal condition of the heat equation is not stable with respect to perturbations in the terminal condition. This issue cannot be resolved without approximation of the inverse with a family of continuous operators, also known as \emph{regularisation}. The research field of \emph{Inverse and Ill-posed Problems} and its branch \emph{Regularisation Theory} focus strongly on the stable approximation of ill-posed and ill-conditioned inverses via \emph{regularisations} \citep{engl1996regularization} and so-called \emph{variational regularisations} \citep{scherzer2009variational,benning2018modern} that are a special class of (nonlinear) regularisations. The optimisation model proposed in \citep{mahendran2015understanding} can be considered as a variational regularisation method with total variation regularisation; however, the work in \citep{mahendran2015understanding} is purely empirical, and to the best of our knowledge no works exist that rigorously prove that the proposed approach is a variational regularisation. 

In this work, we propose a novel regularisation framework based on lifting with tailored Bregman distances and prove that the proposed framework is a convergent variational regularisation for the inverse problem of estimating the inputs from single-layer perceptrons or the inverse problem of estimating hidden variables in a multi-layer perceptron sequentially. While there has been substantial work in previous years that focuses on utilising neural networks as nonlinear operators in variational regularisation methods \citep{lunz2018adversarial,arridge2019solving,Schwab_2019_null_space,li2020nett,mukherjee21}, this is the first work that provides theoretical guarantees for the stable, model-based inversion of neural networks to the best of our knowledge.

Our contributions are three-fold. 1) We propose a novel framework for the regularised inversion of multi-layer perceptrons, respectively feed-forward neural networks, that is based on the lifted Bregman framework recently proposed by the authors in \cite{wang2022lifted}. 2) We show that for the single-layer perceptron case, the proposed variational regularisation approach is a provably convergent regularisation under very mild assumptions. To our knowledge, this is the first time that an inversion method has been proposed that does not just allow to perform inversion empirically, but for which we can prove that the proposed method is a convergent regularisation method without overly restrictive assumptions such as differentiability of the activation function and the presence of a tangential cone condition. 3) We propose a proximal first-order optimisation strategy to solve the proposed variational regularisation method and present several numerical examples that support the effectiveness of the proposed model-based regularisation approach.

The paper is structured as follows. In Section \ref{sec:feed-forward} we introduce the lifted Bregman formulation for the model-based inversion of feed-forward neural networks. In Section \ref{sec:convergence} we prove that for the single-layer perceptron case the proposed model is a convergent variational regularisation method and provide general error estimates as well as error estimates for a concrete example of a perceptron with ReLU activation function. In Section \ref{sec:implementation} we discuss how to implement the proposed variational regularisation computationally for both the single-layer and multi-layer perceptron setting with a generalisation of the primal-dual hybrid gradient method and coordinate descent. Subsequently, we present numerical results that demonstrate empirically that the proposed approach is a model-based regularisation in Section \ref{sec:results}, before we conclude this work with a brief section on conclusions and outlook in Section \ref{sec:outlook}.

\section{Model-based inversion of feed-forward networks}\label{sec:feed-forward}
Suppose we are given an $\layer$-layer feed-forward neural network $\mathcal{N}:\R^n \times \mathcal{P} \rightarrow \R^m$ of the form
\begin{align}
    \mathcal{N}(x,\mathbf{\Theta}) = 
    \sigma_{\layer}(f(\sigma_{\layer-1}(f(\ldots \sigma_{1}(f(x,\Theta_1))\ldots)) ,\Theta_{\layer})), \label{eq:feed-forward-nn}
\end{align}
for input data $x \in \R^n$ and \textbf{pre-trained} parameters $\mathbf{\Theta} \in \mathcal{P}$. Here, $\{ \sigma_{l} \}_{l = 1}^\layer$ denotes the collection of nonlinear activation functions and $f$ denotes a generic function parametrised by parameters $\{ \Theta_l \}_{l = 1}^\layer$. For ease of notation, we use $\mathbf{\Theta}$ to refer to all parameters $\{ \Theta_l \}_{l = 1}^\layer$. For a given network output $y \in \mathbb{R}^m$, our goal is to solve the inverse problem
\begin{align*}
    \mathcal{N}(x,\mathbf{\Theta}) = \data
\end{align*}
for the unknown input $x \in \mathbb{R}^n$. We propose to approximate the inverse of this nonlinear, potentially ill-posed inverse problem via the minimisation of a lifted Bregman formulation of the form
\begin{align}
    \left( \begin{matrix} x^\alpha \\ x_1^\alpha \\ \vdots \\ x_{\layer - 1}^\alpha \end{matrix} \right) \in \argmin_{x, x_1, \ldots, x_{\layer - 1}} \left\{ \sum_{l = 1}^{\layer} \bregmanloss[\Psi_l](x_l, f(x_{l - 1}, \Theta_l)) + \alpha \regfct(x) \right\} \, ,\label{eq:variational-regularisation}
\end{align}
where we assume $x_0 = x$ and $x_\layer = \noisydata$ for simplicity of notation. The data $\noisydata$ is a perturbed version of $\data$, for which we assume $\bregmanloss[\Psi_\layer](\noisydata, f(x_{\layer - 1}^\dagger, \Theta_\layer)) \leq \noisebound^2$, for some constant $\noisebound \geq 0$ and $y = \sigma_{\layer}(f(x_{\layer - 1}^\dagger, \Theta_\layer))$. The functions $\bregmanloss[\Psi_l]$ for $l = 1, \ldots, \layer$ are defined as 
\begin{align}
    \bregmanloss[\Psi_l](x, z) = \frac{1}{2} \| x \|^2 + \Psi_l(x) + \left( \frac12 \| \cdot \|^2 + \Psi_l \right)^\ast(z) - \langle x, z \rangle \, ,\label{eq:bregman-loss}
\end{align}
for a proper, convex and lower semi-continuous function $\Psi_l \colon \mathbb{R}^{n_l} \rightarrow \mathbb{R} \cup \{ \infty \}$. The notation $\left( \frac12 \| \cdot \|^2 + \Psi_l \right)^\ast$ refers to the convex or Fenchel conjugate of $\frac12 \| \cdot \|^2 + \Psi_l$, i.e. $\left( \frac12 \| \cdot \|^2 + \Psi_l \right)^\ast(z) = \sup_y \langle z, y \rangle - \frac12 \| y \|^2 - \Psi_l(y)$. Last but not least, the function $R:\mathbb{R}^n \rightarrow \mathbb{R} \cup \{ \infty \}$ is a proper, convex, and lower semi-continuous function that enables us to incorporate a-priori information into the inversion process. The impact of this is controlled by the parameter $\alpha > 0$.

Please note that the functions $\bregmanloss[\Psi_l]$ have some useful properties and are directly connected to the chosen activation functions $\{ \sigma_l \}_{l = 1}^\layer$. Following \citep{wang2022lifted}, we observe
\begin{align*}
    \bregmanloss[\Psi_l](x, z) \geq \frac12 \| \sigma_l(z) - x \|^2 \, ,
\end{align*}
where $\sigma_l :\mathbb{R}^{n_l} \rightarrow \mathbb{R}^{n_l}$ is the proximal map with respect to $\Psi_l$, i.e.
\begin{align*}
    \sigma_l(z) = \argmin_{y \in \mathbb{R}^{n_l}} \left\{ \frac12 \| y - z \|^2 + \Psi_l(y) \right\} \, ,
\end{align*}
for all $l \in \{1, \ldots, \layer\}$. This means that we will solely focus on feed-forward neural networks with \textbf{nonlinear activation functions that are proximal maps}.

Another useful property is that the functions $\bregmanloss[\Psi_l]$ are continuously differentiable with respect to their second argument. If we define $F_x^l(z) := \bregmanloss[\Psi_l](x, z)$, we observe
\begin{align}
    \nabla F_x^l(z) = \sigma_l(z) - x \, .\label{eq:gradient-property}
\end{align}
Please note that the family of objective functions $\bregmanloss[\Psi_l]$ satisfies several other interesting properties; we refer the interested reader to \citep[Theorem 10]{wang2022lifted}.

For the remainder of this work, we assume that the parametrised functions $f$ are affine-linear in the first argument, with parameters $\Theta_l$. A concrete example is the affine-linear transformation $f(x, \Theta_l) = W_l x + b_l$, for a (weight) matrix $W_l \in \mathbb{R}^{n_l \times n_{l - 1}}$, a (bias) vector $b_l \in \mathbb{R}^{n_l}$ and the collection of parameters $\Theta_l = (W_l, b_l)$. 

In the next section we show that \eqref{eq:variational-regularisation} is a variational regularisation method for $\layer = 1$ and prove a convergence rate with which the solution of \eqref{eq:variational-regularisation} converges towards the true input of a perceptron when $\noisebound$ converges to zero.

\section{Convergence analysis and error estimates}\label{sec:convergence}
In this section we show that the proposed model \eqref{eq:variational-regularisation} is a convergent variational regularisation for the specific choice $\layer = 1$ and the assumption $f(x, \Theta) = Wx + b$ for $\Theta = (W, b)$, which reduces \eqref{eq:variational-regularisation} to a variational regularisation model for the perceptron case studied in \cite{wang2020generalised}. In contrast to \cite{wang2020generalised} we are not interested in estimating the perceptron parameters $W$ and $b$ but assume that these are fixed, and that we study the regularisation operator 
\begin{align}
    \regop\colon \dom(\Psi) \rightrightarrows \R^n \, , \qquad \regop\colon \noisydata \rightrightarrows \reconstruction \in \argmin_{x \in \R^n} \left\{ \bregmanloss\left(\noisydata, W x + b\right) + \alpha \regfct(x) \right\} \, ,\label{eq:perceptron-variational-regularisation}
\end{align}
where $\dom(\Psi)$ is defined as $\dom(\Psi) := \{ y \in \mathbb{R}^m \, | \, \Psi(y) < \infty \}$. We first want to establish under which assumptions \eqref{eq:perceptron-variational-regularisation} is well-defined for all $\noisydata$. 

\subsection{Well-definedness}\label{sec:well-defined}
For simplicity, we focus on the finite-dimensional setting with network inputs in $\mathbb{R}^n$ and outputs in $\dom(\Psi)$. However, the following analysis also extends to more general Banach space settings with additional assumptions on the operator $W$, see for instance \citep[Section 5.1]{benning2018modern}. Following \citep{benning2018modern}, we assume that $\regfct$ is non-negative and the polar of a proper function, i.e. $\regfct = H^\ast$ for a proper function $H:\mathbb{R}^n \rightarrow \mathbb{R} \cup \{ \infty \}$. Note that this automatically implies convexity of $\regfct$. Moreover, we assume that $\Psi$ is a proper, non-negative and convex function that is continuous on $\text{dom}(\Psi)$, which implies that $\bregmanloss$ is proper, non-negative, convex in its second argument and continuous in its first argument for every $\noisydata \in \text{dom}(\Psi)$. Then, for every $g \in \text{dom}(\Psi)$ there exists $x$ with
\begin{align*}
    \bregmanloss(g, Wx + b) + \alpha \regfct(x) < \infty \, .
\end{align*}
Last but not least, we assume that $\regfct$ and $\Psi$ are chosen such that for each $g \in \text{dom}(\Psi)$ and $\alpha > 0$ we have
\begin{align*}
    \| x \| \leq c(a, b, \| g \|), \qquad \text{if} \quad \bregmanloss(g, Wx + b) \leq a \quad \text{and} \quad \alpha \regfct(x) \leq d \, ,
\end{align*}
for constants $a, d$ and a constant $c$ that depends monotonically non-decreasing on all arguments. With these assumptions, we can then verify the following lemma.
\begin{theorem}\label{lem:well-defined}
Let the assumptions outlined in the previous paragraph be satisfied.
\begin{enumerate}
    \item Then, for every $\data \in \left\{ g \in \dom(\Psi) \, \left| \, \argmin_{x, \in \mathbb{R}^n, \regfct(x) < \infty} \bregmanloss(g, Wx + b) \neq \emptyset \right. \right\}$ the selection operator
    \begin{align*}
        \mathcal{S}(\data) = \argmin_{x \in \mathbb{R}^n} \left\{ \regfct(x) \, \left| \, x \in \argmin_{\tilde{x} \in \mathbb{R}^n} \bregmanloss(\data, W\tilde{x} + b) \right. \right\}
    \end{align*}
    is well-defined.\label{lem:selection-operator}
    \item The regularisation operator $\regop$ as defined in \eqref{eq:variational-regularisation} is well-defined in the sense that for every $\data \in \dom(\Psi)$ there exists $x_\alpha \in \mathbb{R}^n$ with $x_\alpha \in \regop(\data)$. Moreover, the set $\regop(\data)$ is a convex set.
    \item For every sequence $\data_n \rightarrow \data \in \dom(\Psi)$ there exists a subsequence $x_{n_k} \in \regop(\data_{n_k})$ converging to an element $x^\ast \in \regop(\data)$.
\end{enumerate}
\begin{proof}
The results follow directly from \citep{benning2018modern}, Lemma 5.5, Theorem 5.6 and Theorem 5.7. The latter statement originally only implies convergence in the weak-star topology; however, since we are in a finite-dimensional Hilbert space, this automatically implies strong convergence here.
\end{proof}
\end{theorem}

\subsection{Error estimates}\label{sec:error}
Having established that \eqref{eq:variational-regularisation} is a regularisation operator, we now want to prove that it is also a convergent regularisation operator in the sense of the estimate
\begin{align}
    \bregmandis(\solution, x^\alpha) \leq C \noisebound \, ,\label{eq:gen-estimate}
\end{align}
such that
\begin{align*}
    \lim_{\delta \rightarrow 0} \sup \left\{ \bregmandis(\solution, x^\alpha) \, \left| \, x^\alpha \in \regop(\noisydata) , \, \noisydata \in \dom(\Psi), \, \bregmanloss(\noisydata, \data) \leq \noisebound^2 \right.\right\} = 0 \, .
\end{align*}
Here, the term $\bregmandis$ denotes the (generalised) Bregman distance (or divergence) (cf. \citep{bregman1967relaxation,kiwiel1997proximal}) with respect to $\regfct$, i.e.
\begin{align*}
    \bregmandis(x, \tilde{x}) = \regfct(x) - \regfct(\tilde{x}) - \langle \tilde{q}, x - \tilde{x} \rangle \, ,
\end{align*}
for two arguments $x, \tilde{x} \in \dom(\regfct)$ and a subgradient $\tilde{q} \in \partial \regfct(\tilde{x}) = \{ q \in \R^n \, | \, \regfct(x) \geq \regfct(\tilde{x}) + \langle q, x - \tilde{x}\rangle , \, \forall \, x \in \dom(\regfct) \}$. The vector $x^\alpha$ is a solution of \eqref{eq:variational-regularisation} with data $\noisydata$ for which we assume $\bregmanloss(\noisydata, \data) \leq \noisebound^2$, and $C \geq 0$ is a constant. The vector $\solution$ is an element of the selection operator as specified in Lemma \ref{lem:well-defined}.\ref{lem:selection-operator}, i.e. $\solution \in \mathcal{S}(\data)$ for $\data \in \dom(\Psi)$. Note that $\solution \in \mathcal{S}(\data)$ is equivalent to $\solution$ being a $\regfct$-minimising vector amongst all vectors that satisfy $0 = W^\ast \left( \actfct(W\solution + b) - \data \right)$, where $\sigma$ denotes the proximal map with respect to $\Psi$. This is due to the fact that $\solution \in \argmin_{\tilde{x} \in \mathbb{R}^n} \bregmanloss(\data, W\tilde{x} + b)$ is equivalent to $0 = \nabla \bregmanloss(\data, W\solution + b) = W^\ast \left( \actfct(W\solution + b) - \data \right)$. Assuming that $\sigma(W\solution + b) - \data$ does not lie in the nullspace of $W^\ast$, this further implies $\data = \actfct(W\solution + b)$.

In order to be able to derive error estimates of the form \eqref{eq:gen-estimate}, we restrict ourselves to solutions $\solution$ that are in the range of $\regop$. This means that there exists $\data^\dagger$ such that $\solution \in \regop(\data^\dagger)$. Considering the optimality condition of \eqref{eq:variational-regularisation} for $\data^\dagger$, this implies
\begin{align*}
    W^\ast \left( \frac{\data^\dagger - \actfct(W\solution + b)}{\alpha} \right) \in \partial \regfct(\solution) \, ,
\end{align*}
which for $\scelem := (\data^\dagger - \actfct(W\solution + b))/\alpha = (\data^\dagger - y)/\alpha$ is equivalent to the existence of a source condition element $\scelem$ that satisfies the source condition (cf. \cite{engl1996regularization,benning2018modern})
\begin{align}
    W^\ast \scelem \in \partial \regfct(\solution) \, , \tag{SC}\label{eq:sc}
\end{align}
In the following, we verify that the symmetric Bregman distance with respect to $\regfct$ between a solution of the regularisation operator and the solution of the inverse problem is converging to zero if the error in the data is converging to zero. The symmetric Bregman distance or Jeffreys distance between two vectors $x$ and $\tilde{x}$ simply is the sum of two Bregman distances with interchanged arguments, i.e.
\begin{align*}
    \symbreg(x, \tilde{x}) := \bregmandis(x, \tilde{x}) + \bregmandis(\tilde{x}, x) = \langle x - \tilde{x}, q - \tilde{q} \rangle \, ,
\end{align*}
for $q \in \partial \regfct(x)$ and $\tilde{q} \in \partial \regfct(\tilde{x})$; hence, an error estimate in the symmetric Bregman distance also implies an error estimate in the classical Bregman distance.

Before we begin our analysis, we recall the concept of the \emph{Jensen-Shannon divergence} \citep{lin1991divergence}, which for general proper, convex and lower semi-continuous functions $F:\R^n \rightarrow \R \cup \{ \infty \}$ generalises to so-called \emph{Burbea-Rao divergences} \citep{burbea1982convexity1,burbea1982convexity2,nielsen2011burbea} and are defined as follows. 
\begin{definition}[Burbea-Rao divergence]
Suppose $F:\R^n \rightarrow \R \cup \{ \infty \}$ is a proper, convex and lower semi-continuous function. The corresponding Burbea-Rao divergence is defined as
\begin{align}
    \jensen[F](x, \tilde{x}) := \frac{1}{2} \left( F(x) + F(\tilde{x}) - 2 F\left( \frac{x + \tilde{x}}{2} \right) \right) \, ,\label{eq:jensen-shannon}
\end{align}
for all $x, \tilde{x} \in \dom(F)$.
\end{definition}
Another important concept that we need in order to establish error estimates is that of Fenchel conjugates (cf. \cite{beck2017first}).
\begin{definition}[Fenchel conjugate]
The Fenchel (or convex) conjugate $F^\ast:\R^n \rightarrow \R \cup \{ -\infty, \infty\}$ of a function $F:\R^n \rightarrow \R \cup \{ -\infty, +\infty\}$ is defined as
\begin{align*}
    F^\ast(w) \colon = \sup_{x \in \R^n} \langle x, w \rangle - F(x) \, .
\end{align*}
\end{definition}
The Fenchel conjugate that is of particular interest to us is the conjugate of the function $\bregmanloss(y, z)$ with respect to the second argument, which we characterise with the following lemma.
\begin{lemma}\label{lem:conjugate}
The Fenchel conjugate of $F_y(z) := \bregmanloss(y, z)$ with respect to the second argument $z$ reads
\begin{align*}
    F_y^\ast(w) = \left( \frac12 \| \cdot \|^2 + \Psi\right)(y + w) - \left( \frac12 \| \cdot \|^2 + \Psi\right)(y) \, .
\end{align*}
\begin{proof}
    From the definition of the Fenchel conjugate we observe
    \begin{align*}
        F_y^\ast(w) &= \sup_{z \in \R^m} \, \langle z, w \rangle - F_y(z) \\
        &= \sup_{z \in \R^m} \, \langle z, w \rangle - \left( \frac12 \| \cdot \|^2 + \Psi\right)(y) - \left( \frac12 \| \cdot \|^2 + \Psi\right)^\ast(z) + \langle y, z \rangle\\
        &= - \left( \frac12 \| \cdot \|^2 + \Psi\right)(y) + \sup_{z \in \R^m} \, \langle z, w + y \rangle - \left( \frac12 \| \cdot \|^2 + \Psi\right)^\ast(z) \\
        &= - \left( \frac12 \| \cdot \|^2 + \Psi\right)(y) + \left( \frac12 \| \cdot \|^2 + \Psi\right)(w + y) \, ,
    \end{align*}
    which concludes the proof.
\end{proof}
\end{lemma}
Having defined the Burbea-Rao divergence and having established the Fenchel conjugate of $\bregmanloss(y, z)$ with respect to the second argument $z$ for fixed $y$, we can now present and verify our main result that is motivated by \citep{benning2011error}.
\begin{theorem}
Suppose $\regfct$ and $\Psi$ satisfy the assumptions outlined in Section \ref{sec:well-defined}. Then, for data $\noisydata$ and $\solution$ that satisfy $\bregmanloss(\noisydata, W\solution + b) \leq \noisebound^2$ with $\noisebound \geq 0$, a solution $x^\alpha \in \regop(\noisydata)$ of the variational regularisation problem \eqref{eq:variational-regularisation}, and a solution $\solution$ of the perceptron problem $\data = \actfct(W\solution + b)$ that satisfies $\solution \in \mathcal{S}(\data)$ and \eqref{eq:sc}, we observe the error estimate
\begin{align}
    \begin{split}
    (1 - c) \bregmanloss(\noisydata, W\reconstruction + b) + \alpha \symbreg(\reconstruction, \solution) {} \leq {} &(1 + c) \noisebound^2 + \frac{\alpha^2}{c} \| \scelem \|^2 \\
    &+ 2c \jensen\left(\noisydata + \frac{\alpha}{c} \scelem, \noisydata - \frac{\alpha}{c} \scelem \right)
    \end{split} \, ,\label{eq:perceptron-error-estimate}
\end{align}
for a constant $c \in (0, 1]$.
\begin{proof}
Every solution $\reconstruction$ that satisfies $\reconstruction \in \regop(\noisydata)$ can equivalently be characterised by the optimality condition
\begin{align*}
    W^\ast \left( \actfct( W \reconstruction + b) - \noisydata \right) + \alpha p_\alpha = 0 \, ,
\end{align*}
for any subgradient $p_\alpha \in \partial \regfct(\reconstruction)$. Subtracting $p^\dagger \in \partial \regfct(\solution)$ from both sides of the equation and taking a dual product with $\reconstruction - \solution$ then yields
\begin{align}
    \langle \actfct( W \reconstruction + b) - \noisydata, W\reconstruction - W\solution \rangle + \alpha \symbreg(\reconstruction, \solution) = - \alpha \langle p^\dagger, \reconstruction - \solution \rangle \, .\label{eq:main-estimate-intermediate-step1}
\end{align}
We easily verify 
\begin{align*}
    \bregmandis[\bregmanloss(\noisydata, W \cdot + b)](\solution, \reconstruction) = \bregmanloss(\noisydata, W\solution + b) - \bregmanloss(\noisydata, W\reconstruction + b) - \langle \actfct( W \reconstruction + b) - \noisydata, W\solution - W\reconstruction \rangle \, ;
\end{align*}
hence, we can replace $\langle \actfct( W \reconstruction + b) - \noisydata, W\reconstruction - W\solution \rangle$ with $\bregmandis[\bregmanloss(\noisydata, W \cdot + b)](\solution, \reconstruction) + \bregmanloss(\noisydata, W\reconstruction + b) - \bregmanloss(\noisydata, W\solution + b)$ in \eqref{eq:main-estimate-intermediate-step1} to obtain 
\begin{align*}
    \bregmandis[\bregmanloss(\noisydata, W \cdot + b)](\solution, \reconstruction) + \bregmanloss(\noisydata, W\reconstruction + b) + \alpha \symbreg(\reconstruction, \solution) = \bregmanloss(\noisydata, W\solution + b) - \alpha \langle p^\dagger, \reconstruction - \solution \rangle \, .
\end{align*}
We know $0 \leq \bregmandis[\bregmanloss(\noisydata, W \cdot + b)](\solution, \reconstruction)$ due to the convexity of $\bregmanloss(\noisydata, W \cdot + b)$, and we also know that \eqref{eq:sc} enables us to choose $p^\dagger = W^\ast \scelem$. Hence, we can estimate
\begin{align*}
    \bregmanloss(\noisydata, W\reconstruction + b) + \alpha \symbreg(\reconstruction, \solution) \leq \bregmanloss(\noisydata, W\solution + b) - \alpha \langle \scelem, W\reconstruction - W\solution \rangle \, .
\end{align*}
Next, we introduce the constant $c \in (0, 1]$ to split the loss functions $\bregmanloss(\noisydata, W\reconstruction + b)$ and $\bregmanloss(\noisydata, W\solution + b)$ into $(1 - c)\bregmanloss(\noisydata, W\reconstruction + b) + c \bregmanloss(\noisydata, W\reconstruction + b)$ and $(1 + c) \bregmanloss(\noisydata, W\solution + b) - c \bregmanloss(\noisydata, W\solution + b)$, respectively. This means we estimate
\begin{align*}
    (1 - c) \bregmanloss(\noisydata, W\reconstruction + b) + \alpha \symbreg(\reconstruction, \solution) {} \leq {} &(1 + c) \bregmanloss(\noisydata, W\solution + b) \\
    &+ \langle \alpha \scelem, W\solution + b \rangle - c \bregmanloss(\noisydata, W\solution + b) \\
    &- \langle \alpha \scelem, W\reconstruction + b \rangle - c \bregmanloss(\noisydata, W\reconstruction + b) \, .
\end{align*}
Next, we make use of Lemma \ref{lem:conjugate} to estimate
\begin{align*}
    \langle \alpha \scelem, W\solution + b \rangle - c \bregmanloss(\noisydata, W\solution + b) &\leq c \left( \left(\frac12 \| \cdot \|^2 + \Psi\right)\left(\noisydata + \frac{\alpha}{c} \scelem \right) - \left(\frac12 \| \cdot \|^2 + \Psi\right)\left(\noisydata\right) \right) \, ,
\intertext{and}
    - \langle \alpha \scelem, W\reconstruction + b \rangle - c \bregmanloss(\noisydata, W\reconstruction + b) &\leq c \left( \left(\frac12 \| \cdot \|^2 + \Psi\right)\left(\noisydata - \frac{\alpha}{c} \scelem \right) - \left(\frac12 \| \cdot \|^2 + \Psi\right)\left(\noisydata\right) \right) \, .
\end{align*}
Adding both estimates together yields
\begin{align*}
    &\langle \alpha \scelem, W\solution + b \rangle - c \bregmanloss(\noisydata, W\solution + b) - \langle \alpha \scelem, W\reconstruction + b \rangle - c \bregmanloss(\noisydata, W\reconstruction + b) \, , \\
    {} \leq {} &\frac{\alpha^2}{c}\| \scelem \|^2 + c \left( \Psi\left(\noisydata + \frac{\alpha}{c} \scelem \right) + \Psi\left(\noisydata - \frac{\alpha}{c} \scelem \right) - 2 \Psi(\noisydata) \right) \, ,\\
    {} = {} & \frac{\alpha^2}{c}\| \scelem \|^2 + 2c \jensen\left(\noisydata + \frac{\alpha}{c} \scelem, \noisydata - \frac{\alpha}{c} \scelem \right) \, ,
\end{align*}
which together with the error bound $\bregmanloss(\noisydata, W\solution + b) \leq \noisebound^2$ concludes the proof.
\end{proof}
\end{theorem}
\begin{remark}
    We want to emphasise that for continuous $\Psi$ and $c > 0$ we automatically observe
    \begin{align*}
        \lim_{\alpha \rightarrow 0} \jensen\left(\noisydata + \frac{\alpha}{c} \scelem, \noisydata - \frac{\alpha}{c} \scelem \right) = 0 \, ,
    \end{align*}
    in which case the important question from an error estimate point-of-view is if the term converges quicker to zero than $\alpha$, as we would need to guarantee $\lim_{\alpha \rightarrow 0} \jensen\left(\noisydata + \frac{\alpha}{c} \scelem, \noisydata - \frac{\alpha}{c} \scelem \right)/\alpha = 0$ in order to guarantee that the symmetric Bregman distance in \eqref{eq:perceptron-error-estimate} converges to zero for $\alpha \rightarrow 0$.
\end{remark}
\begin{example}[ReLU perceptron]
    Let us consider a concrete example to demonstrate that \eqref{eq:perceptron-variational-regularisation} is a convergent regularisation with respect to the symmetric Bregman distance of $R$. We know that for $\actfct(z) = \prox(z) = \max(0, z)$ to hold true we have to choose $\Psi(z) = \begin{cases} 0 & z \in [0, \infty)^m \\ \infty & \text{else}\end{cases}$. This means that for $\bregmanloss(\noisydata, z)$ to be well-defined for any $z$ we require $\noisydata_i \geq 0$ for all $i \in \{1, \ldots, m\}$. In order for the Burbea-Rao divergence to be well-defined, we further require 
    \begin{align*}
        -\frac{c}{\alpha} \noisydata_i \leq \scelem_i \leq \frac{c}{\alpha} \noisydata_i \, ,
    \end{align*}
    for all $i \in \{1, \ldots, m\}$, or $\| \scelem \|_\infty \leq (c \| \noisydata \|_\infty / \alpha)$ in more compact notation. If $\| \scelem \|_\infty \leq (c \| \noisydata \|_\infty / \alpha)$ is guaranteed, we observe $\jensen\left(\noisydata + \frac{\alpha}{c} \scelem, \noisydata - \frac{\alpha}{c} \scelem \right) = 0$. Hence, we can simplify the estimate \eqref{eq:perceptron-error-estimate} to
    \begin{align*}
        \frac{1 - c}{\alpha} \bregmanloss(\noisydata, W\reconstruction + b) + \symbreg(\reconstruction, \solution) \leq \frac{1 + c}{\alpha} \noisebound^2 + \frac{\alpha}{c} \| \scelem \|^2 \, ,
    \end{align*}
    where we have also divided by $\alpha$ on both sides of the inequality. If we choose $\alpha(\noisebound) = \sqrt{c(1 + c)}\noisebound/\| \scelem \|$, we obtain the estimate
    \begin{align*}
        \frac{(1 - c)\|\scelem\|}{\delta \sqrt{c (1 + c)}} \bregmanloss(\noisydata, W\reconstruction[\alpha(\noisebound)] + b) + \symbreg(\reconstruction[\alpha(\noisebound)], \solution) \leq 2 \sqrt{\frac{1 + c}{c}} \| \scelem \| \noisebound \, ,
    \end{align*}
    as long as we can ensure 
    \begin{align*}
        \left\| \frac{\scelem}{\|\scelem \|} \right\|_\infty \leq \sqrt{\frac{c}{1 + c}} \left\| \frac{\noisydata}{\noisebound} \right\|_\infty \, .
    \end{align*}
    Together with $\symbreg(\reconstruction[\alpha(\noisebound)], \solution) \geq \bregmandis(\solution, \reconstruction[\alpha(\noisebound)])$ we have established an estimate of the form \eqref{eq:gen-estimate}, with constant $C = 2 \sqrt{\frac{1 + c}{c}} \| \scelem \|$. Hence, we have verified that the variational regularisation method \eqref{eq:variational-regularisation} is not only a regularisation method but even a convergent regularisation method in this specific example.
\end{example}
We want to briefly comment on the extension of the convergence analysis to the general case $\layer > 1$ with the following remark. 
\begin{remark}
    The presented convergence analysis easily extends to a sequential, layer-wise inversion approach. Suppose we have $\layer$ layers and begin with the final layer, then we can formulate the variational problem
    \begin{align*}
        x_{\layer - 1}^\alpha &\in \argmin_{x_{\layer - 1}} \left\{ \bregmanloss[\Psi_\layer](\noisydata, W_{\layer} x_{\layer - 1} + b_{\layer}) + \alpha_{\layer - 1} \Psi_{\layer - 1}(x_{\layer - 1}) \right\} \, ,
    \end{align*}
    which is also of the form of \eqref{eq:perceptron-variational-regularisation}, but where $\regfct$ has been replaced with $\Psi_{\layer - 1}$. Alternatively, one can also replace $\Psi_{\layer - 1}$ with another function $\regfct_{\layer - 1}$ if good prior knowledge for the auxiliary variable $x_{\layer - 1}$ exists. Once we have estimated $x_{\layer - 1}^\alpha$, we can recursively estimate
    \begin{align*}
        x_{l - 1}^\alpha &\in \argmin_{x_{l - 1}} \left\{ \bregmanloss[\Psi_l](x_l^\alpha, W_{l} x_{l - 1} + b_{l}) + \alpha_{l - 1} \Psi_{l - 1}(x_{l - 1}) \right\} \, ,
    \end{align*}
    for $l = \layer - 1, \ldots, 2$ and subsequently compute $\reconstruction$ as a solution of \eqref{eq:variational-regularisation} but with data $x_1^\alpha$ instead of $\noisydata$. 
    
    The advantage of such a sequential approach is that every individual regularisation problem is convex and the previously presented theorems and guarantees still apply. The disadvantage is that for this approach to work in theory, we require bounds for every auxiliary variable of the form $\bregmanloss[\Psi_l](x_l^\alpha, W_l x_{l - 1}^\alpha + b_l) \leq \delta_l^2$, which is a rather unrealistic assumption. Moreover, it is also not realistic to assume that good prior knowledge for the auxiliary variables exist. 
    
    Please note that showing that the simultaneous approach \eqref{eq:variational-regularisation} is a (convergent) variational regularisation is beyond the scope of this work as it is harder and potentially requires additional assumptions for the following reason. The overall objective function in \eqref{eq:variational-regularisation} is no longer guaranteed to be convex with respect to all variables simultaneously, which means that we cannot simply carry over the analysis of the single-layer to the multi-layer perceptron case.
    
\end{remark}

This concludes the theoretical analysis of the perceptron inversion model. In the following section we focus on how to implement \eqref{eq:perceptron-variational-regularisation} and its more general counterpart \eqref{eq:variational-regularisation}.





\section{Implementation}\label{sec:implementation}

In this section, we describe how to computationally implement the proposed variational regularisation for both the single-layer and the multi-layer perceptron setting. More specifically, we show that the proposed variational regularisation can be efficiently solved via a generalised primal-dual hybrid gradient method and a coordinate descent approach. 


\subsection{Inverting perceptrons}\label{sec:implementation-perceptron}

To begin with, we first consider the example of inverting a (single-layer) perceptron. For $\layer=1$, Problem \eqref{eq:variational-regularisation} reduces to \eqref{eq:perceptron-variational-regularisation}, which for a composite regularisation function $\regfct \circ K$ reads
\begin{align}
    x^\alpha \in \argmin_x \left\{ \bregmanloss\left(\noisydata, f(x,\Theta) \right) + \alpha \regfct(Kx) \right\} \, . \label{eq:generalised-perceptron-inversion}
\end{align}
Here $K$ is a matrix and $\alpha \regfct(Kx)$ denotes the regularisation function acting on the argument $x$. The above Problem \eqref{eq:generalised-perceptron-inversion} can be reformulated to the saddle-point problem
\begin{align}
    \min_x \max_z \bregmanloss\left(\noisydata, f(x,\Theta)\right) + \langle z, Kx \rangle - \alpha \regfct^\ast(z) \, , \label{eq:perceptron-saddle-point-problem}
\end{align}
where $\regfct^\ast$ denotes the convex conjugate of $\regfct$. Computationally, we can then solve the saddle-point problem with a generalised primal-dual hybrid gradient (PDHG) method \citep{zhu2008efficient,pock2009algorithm,esser2010general,chambolle2011first,chambolle2016introduction,benning2021bregman}:
\begin{subequations}\label{eq:perceptron-inversion-pdhg}
    \begin{align}
    x^{k + 1} &= x^k - \tau_{x} \left( \left( \text{prox}_{\Psi} \left(f(x^k, \Theta) \right) - \noisydata \right) \mathcal{J}_f^x(x^k, \Theta) + \alpha K^\top z^k \right) \, \label{eq:perceptron-inversion-pdhg-1}\\
    z^{k + 1} &= \text{prox}_{\tau_z R^\ast}\left( z^k + \tau_z \alpha K  \left(  2 x^{k + 1} -  x^{k} \right) \right) \, . 
\end{align}
\end{subequations}
where we alternate between a descent step in the $x$ variable and an ascent step in the dual variable $z$. Since \eqref{eq:generalised-perceptron-inversion} is a convex minimisation problem, \eqref{eq:perceptron-inversion-pdhg} is guaranteed to converge globally for arbitrary starting point, given that $\tau_x$ and $\tau_z$ are chosen such that $\tau_x \tau_z < 1/\| K \|$.

In this work, we will focus on the discrete total variation $\| \nabla x \|_{p, 1}$, \citep{rudin1992nonlinear,chambolle1997image}, as our regularisation function $\regfct(Kx)$, but other choices are certainly possible. If we consider a two-dimensional scalar-valued image $x\in \R^{H \times W}$, we can define a finite forward difference discretisation of the gradient operator $\nabla: \R^{H \times W} \rightarrow \R^{H \times W \times 2}$ as
\begin{align*}
    &(\nabla x)_{i,j,1} = \begin{cases}
        x_{i+1,j} - x_{i,j} \;\; \text{if } 1 \leq i < H, \\
        0 \;\; \text{else,}\\
    \end{cases} \, ,
    &(\nabla x)_{i,j,2} = 
    \begin{cases}
        x_{i,j+1} - x_{i,j} \;\; \text{if } 1 \leq j < W, \\
        0 \;\; \text{else.}
    \end{cases} \, .
\end{align*}
The discrete total variation is defined as the $\ell_1$ norm of the $p$-norm of the pixel-wise image gradients, i.e. 
\begin{align*}
    \|\nabla x \|_{p,1} = \sum_{i=1}^{H}\sum_{j=1}^{W} |\left(\nabla x\right)_{i,j}|_{p} = \sum_{i=1}^{H}\sum_{j=1}^{W} 
    \left(\left(\nabla x\right)_{i,j,1}^{p} + \left(\nabla x\right)_{i,j,2}^{p}\right)^{1/p} \, .
\end{align*}%
For our numerical results we consider the isotropic total variation and consequently choose $p=2$. Hence for a perceptron with affine-linear transformation $f(x, \Theta) = Wx + b$, and with $\sigma = \prox$ denoting the activation function, the PDHG approach \eqref{eq:perceptron-inversion-pdhg} of solving the perceptron inversion problem \eqref{eq:perceptron-variational-regularisation} can be summarised as
\begin{subequations}\label{eq:perceptron-inversion-tv-pdhg}
    \begin{align}
    x^{k + 1} &= x^k - \tau_x \left( W^\top \left( \sigma(Wx^k + b) - y \right) - \alpha \text{div} z^k \right) \, , \label{eq:perceptron-inversion-tv-pdhg-1}\\
    z^{k + 1} &= \text{prox}_{\tau_z \| \cdot \|_{2, 1}^\ast}\left( z^k + \tau_z \left( 2 \alpha \nabla x^{k + 1} - \alpha \nabla x^{k} \right) \right) \, . \label{eq:perceptron-inversion-tv-pdhg-2}
\end{align}
\end{subequations}
Please note that we define the discrete approximation of the divergence $\text{div}$ such that it satisfies $\text{div} = - \nabla^\top$ in order to be the negative transpose of the discretised finite difference approximation of the gradient in analogy to the continuous case, which is why the sign in \eqref{eq:perceptron-inversion-tv-pdhg-1} is flipped in comparison to \eqref{eq:perceptron-inversion-pdhg-1}. The proximal map with regards to the convex conjugate of $\| \cdot \|_{2, 1}^\ast$ is simply the argument itself if the maximum of the Euclidean vector-norm per pixel is bounded by one or the projection onto this unit ball.


\subsection{Inverting multi-layer perceptrons}\label{sec:implementation-mlp}

We now discuss the implementation of the inversion of multi-layer perceptrons with $\layer$ layers as described in \eqref{eq:variational-regularisation}. Note that in this case in order to minimise for $x$, we also need to optimise with respect to the auxiliary variables $x_1, \ldots, x_{\layer-1}$. 

For the minimisation of \eqref{eq:variational-regularisation} we consider an alternating minimisation approach, also known as \emph{coordinate descent} \citep{beck2013convergence,wright2015coordinate,wright2022optimization}. In this approach we minimise the objective with respect to one variable at a time. In particular, we focus on a semi-explicit coordinate descent algorithm, where we linearise with respect to the smooth functions of the overall objective function. 
This breaks down the overall minimisation problem into $\layer$ sub-problems, where for $x_0$ and each $x_l$ variable for $l \in \{1, \ldots, \layer-1\}$, we have individual minimisation problems of the following form:
\begin{subequations}
\begin{align}
    x_0^{k + 1} {} = {} &\argmin_{x_0} \left\{ \left( \frac12 \| \cdot \|^2 + \Psi_{1}\right)^\ast\left(f(x_0, \Theta_{1}) \right) - \left\langle x_1^{k}, f(x_{0}, \Theta_1) \right\rangle  + \alpha R(K x_0) \right\} 
    \, , \label{eq:coordinate-descent-first-variable} \\
    \begin{split}
    x_l^{k + 1} {} = {} &\argmin_{x_l} \left\{ \left( \frac12 \| \cdot \|^2 + \Psi_{l}\right)(x_l) - \left\langle x_l, f(x_{l - 1}^{k + 1}, \Theta_l) \right\rangle + \frac{1}{2\tau_{x_l}} \| x_l - x_l^k \|^2 \right.\\
    &\qquad \qquad \left. + \left\langle x_l, \left( \text{prox}_{\Psi_{l + 1}} \left(f(x_l^k, \Theta_{1 + 1}) \right) - x_{l + 1}^k \right) \mathcal{J}_f^x(x_l^k, \Theta_{1 + 1}) \right\rangle \vphantom{\left( \frac12 \| \cdot \|^2 + \Psi_{l}\right)(x_l)} \right\} \, .
    \end{split}\label{eq:coordinate-descent-all-others}
\end{align}\label{eq:coordinate-descent-all-variables}%
\end{subequations}
Note that one advantage for adopting this approach is that we exploit that the overall objective function is convex in each individual variable when all other variables are kept fixed. In the following, we will discuss different strategies to computationally solve each sub-problem. 

When optimising with respect to the input variable $x_0$, the structure of sub-problem \eqref{eq:coordinate-descent-first-variable} is identical to the perceptron inversion problem that we have discussed in Section \ref{sec:implementation-perceptron}. Hence, we can approximate $x_0^{k + 1}$ with \eqref{eq:perceptron-saddle-point-problem}, but now with respect to $x_1^k$ instead of $\noisydata$, which yields the iteration
\begin{subequations}
  \begin{align}
    x_0^{k + 1} &= x_0^k - \tau_{x_0} \left( \left( \text{prox}_{\Psi} \left(f(x_0^k, \Theta_{1}) \right) - x_{1}^k \right) \mathcal{J}_f^x(x_0^k, \Theta_{1}) + \alpha K^\top z^k \right) \, , \\
    z^{k + 1} &= \text{prox}_{\tau_z R^\ast}\left( z^k + \tau_z \alpha K  \left(  2 x^{k + 1} -  x^{k} \right) \right) \, .
\end{align}\label{eq:mlp-x0-pdhg}%
\end{subequations}
For each auxiliary variable $x_{l}$ with $l \in \{1, \ldots, \layer-1\}$, the sub-problem associated with \eqref{eq:coordinate-descent-all-others} amounts to solving a proximal gradient step with suitable step-size $\tau_{x_l}$, which we can rewrite to 
\begin{align}\label{eq:x_l-proximal-gradient-specific-form}
\begin{split}
    x_l^{k + 1} {} = {} &\text{prox}_{\frac{\tau_{x_l}}{1 + \tau_{x_l}} \Psi_l}\left( \frac{1}{1 + \tau_{x_l}} \left( x_l^k - \tau_{x_l} \left( \left( \text{prox}_{\Psi_l}\left( f(x_l^k, \Theta_{l + 1}) \right) - x_{l + 1}^k \right) \mathcal{J}_f^x(x_l^k, \Theta_{l + 1}) \right. \right. \right. \\
    &\qquad \qquad \qquad \left. \left. \left. - f(x_{l - 1}^{k + 1}, \Theta_l) \right) \right) \right) 
\end{split} \, .
\end{align}

This concludes the discussion on the implementation of the regularised single-layer and multi-layer perceptron inversion. In the next section, we present some numerical results to demonstrate the effectiveness of the proposed approaches empirically. 


\section{Numerical results}\label{sec:results}

In this section, we present numerical results for the perceptron inversion problem implemented with the PDHG algorithm as outlined in \eqref{eq:perceptron-inversion-tv-pdhg}, and for the multi-layer perceptron inversion problem implemented with the coordinate descent approach as described in \eqref{eq:mlp-x0-pdhg} and \eqref{eq:x_l-proximal-gradient-specific-form}. All results have been computed using PyTorch 3.7 on an Intel Xeon CPU E5-2630 v4.

\subsection{The Perceptron}
\label{sec:perceptron-numerical-results}
We present results for two experiments: the first one is the perceptron inversion of the image of a circle from the noisy output of the perceptron, where we compare the Landweber regularisation and the total-variation-based variational regularisation \eqref{eq:perceptron-variational-regularisation}. For the second experiment, we perform perceptron inversion for samples from the MNIST dataset and compare them with the performance of linear and nonlinear decoders.   

\begin{figure}[!t]
    \centering
    \begin{minipage}{0.325\textwidth}
        \centering 
        \includegraphics[width = \linewidth]{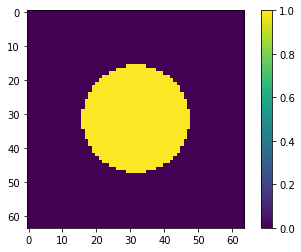}
        \caption{Groundtruth image $\solution$ of a circle.}
        \label{fig:toy-example-circle-groundtruth}
    \end{minipage}
    \hspace{0.0001\textwidth}
    \begin{minipage}{0.333\textwidth}
        \centering 
        \includegraphics[width = \linewidth]{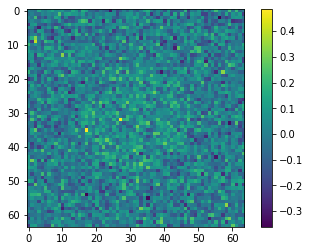}
        \caption{Inverted image via Landweber regularisation.}
        \label{fig:toy-example-circle-landweber}
    \end{minipage}
    \begin{minipage}{0.315\textwidth}
        \centering 
        \includegraphics[width = \linewidth]{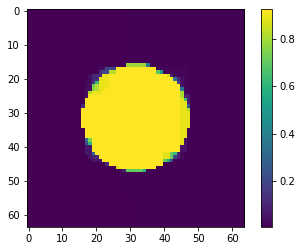}
        \caption{Inverted image via TV regularisation.}\label{fig:toy-example-circle-tv}
    \end{minipage}    
\end{figure}

\noindent \textbf{Circle} We begin with the toy example of recovering the image of a circle from noisy measurements of a ReLU perceptron. To prepare the experiment, we generate a circle image $\solution \in \R^{64\times 64}$, as shown in Figure \ref{fig:toy-example-circle-groundtruth}. We construct a perceptron with ReLU activation function using random weights and biases where $W \in \R^{512 \times 4096}, b \in \R^{512 \times 1}$. The weights operates on the column-vector representation of x, where $x \in \R^{4096 \times 1} $. The noise-free data is generated via the forward operation of the model, i.e. $\data = \sigma(W\solution+b)$. We generate noisy data $\noisydata$ by adding Gaussian noise with mean 0 and standard deviation $0.005$. Note that we clip all the negative values of $\noisydata$ to ensure $\noisydata \in \dom(\Psi)$.

A first attempt to solve this ill-posed perceptron inversion problem is via Landweber regularisation \citep{landweber1951iteration}. In Figure \ref{fig:toy-example-circle-landweber} we see the reconstructed image obtained with Landweber regularisation in combination with early stopping following Morozov's discrepancy principle \cite{morozov2012methods,engl1996regularization}. Even though the Landweber regularised reconstruction matches the data up to the discrepancy value $\| \actfct(Wx^K + b) - \noisydata\|$, the recovered image does not resemble the image $\solution$. We will discuss shortly the reason for this visually poor inversion. In comparison, we see a regularised inversion via the total variation regularisation approach following \eqref{eq:perceptron-inversion-tv-pdhg} in Figure \ref{fig:toy-example-circle-tv}. The regularisation parameter for this reconstruction is chosen as $\alpha=1.5 \times 10^{-2}$.
Both $x_0$ and $z$ are initialised with zero vectors. The stepsize-parameters are chosen as $\tau_{x} = 1.99/\|W\|^2_2$ and $\tau_z = 1/(8\alpha)$, see \citep{chambolle2004algorithm}. We stop the iterations when changes in $x_0$ and $z$ in norm are less than a threshold of $10^{-5}$ or when we reach the maximum number of iterations, which we set to $10000$. As shown in Figure \ref{fig:toy-example-circle-tv}, the TV-regularisation approach is capable of finding a (visually) more meaningful solution. 

To explain why the Landweber iteration performs worse compared to the total variation regularisation for this specific example, we compare the $\ell_2$ norms of each two solutions and the groundtruth image $\solution$. The $\ell_2$ norm of the Landweber solution in Figure \ref{fig:toy-example-circle-landweber} measures $6.58$ while the TV-regularised solution as in Figure \ref{fig:toy-example-circle-tv} and the groundtruth image $\solution$ measure $25.69$ and $28.07$ respectively. This is not surprising, as the Landweber iteration is known to converge to a minimal Euclidean norm solution if the noise level converges to zero. On the other hand, when we compare the TV semi-norm of each solution, the groundtruth image in measures $128.0$, while the Landweber solution in Figure \ref{fig:toy-example-circle-landweber} and TV-regularised solution in Figure \ref{fig:toy-example-circle-tv} measure $707.02$ and $114.93$ respectively, suggesting that the TV-semi-norm is a more suitable regularisation function for the inversion of cartoon-like images such as $\solution$. 

\begin{figure}[!t]
\centering
\begin{minipage}{0.45\textwidth}
        \centering 
        \includegraphics[width = \textwidth]{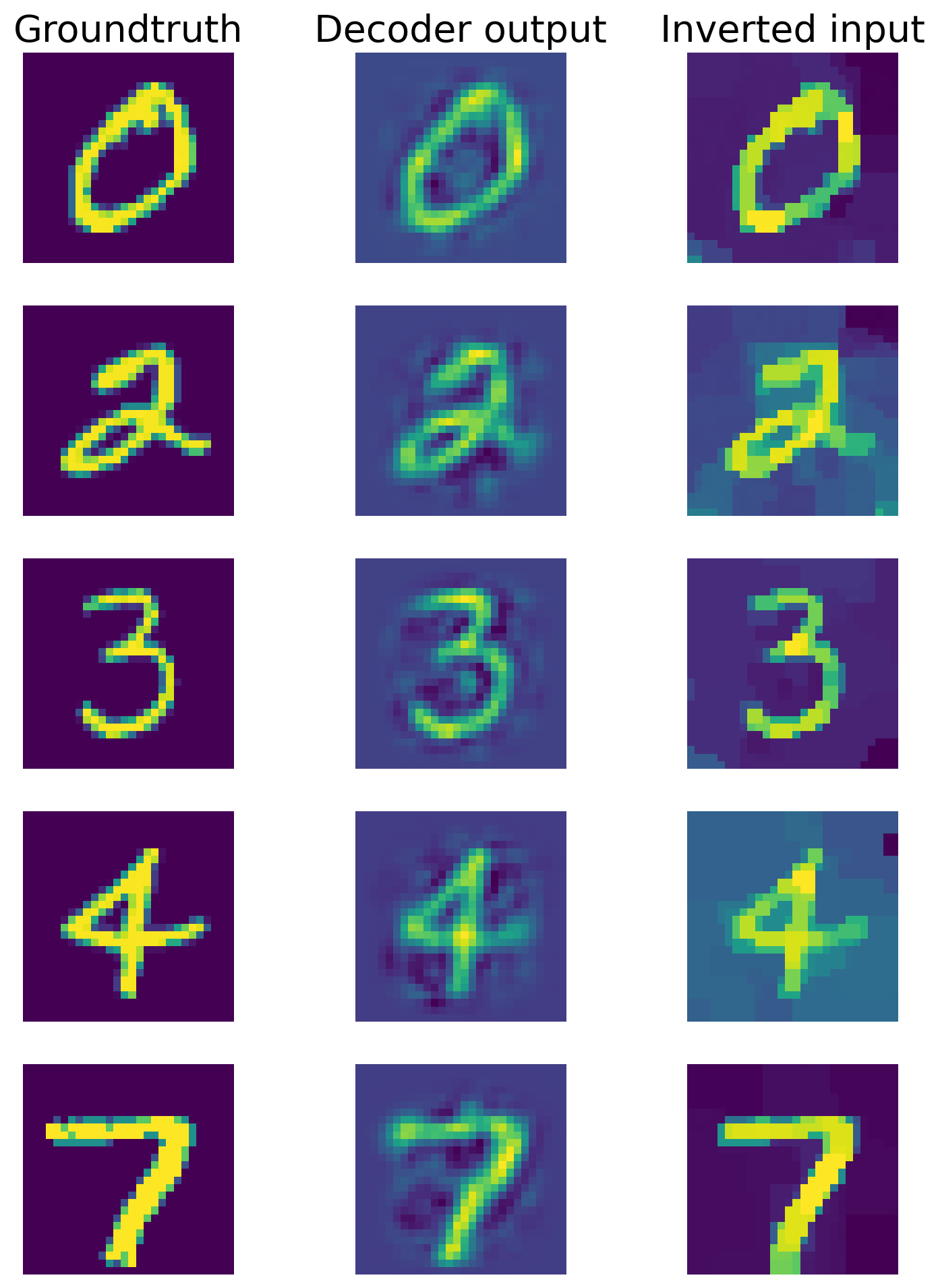}
        \caption{Groundtruth input images from the MNIST training dataset, together with the corresponding autoencoder output images and inverted input images of the perceptron}
        \label{fig:perceptron-inversion-visualisation-train-images}
\end{minipage}
\begin{minipage}{0.45\textwidth}
        \centering
        \includegraphics[width = \textwidth]{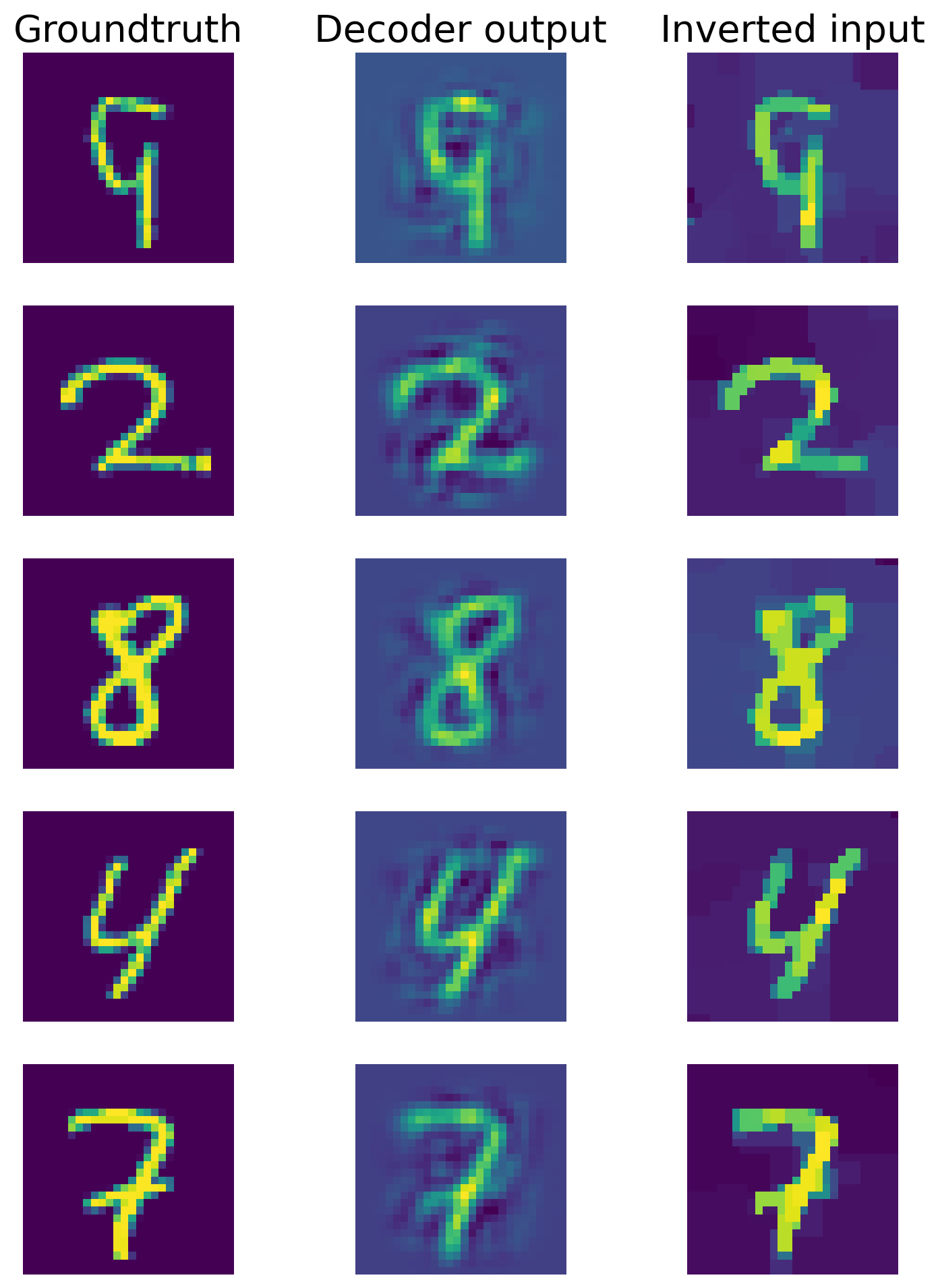}
        \caption{Groundtruth input images from the MNIST validation dataset, together with the corresponding autoencoder output images and inverted input images of the perceptron}
        \label{fig:perceptron-inversion-visualisation-val-images}
\end{minipage}
\end{figure}

\noindent \textbf{MNIST} In this second example, we perform perceptron inversion on the MNIST dataset \citep{lecun1998gradient}. In particular, we consider the following experimental setup. We first train an autoencoder $\mathcal{A}(x) = \mathcal{D}(\mathcal{E}(x,\mathbf{\Theta_{\mathcal{E}}}),\mathbf{\Theta_{\mathcal{D}}})$, where $\mathcal{D}(\cdot,\mathbf{\Theta}_{\mathcal{D}})$ and $\mathcal{E}(\cdot,\mathbf{\Theta_{\mathcal{E}}})$ denotes the decoder and the encoder, parametrised by parameters $\mathbf{\Theta}_{\mathcal{D}}$ and $\mathbf{\Theta}_{\mathcal{E}}$ respectively. We pre-train the autoencoder $\mathcal{A}$, compute the code $\mathcal{E}(x,\mathbf{\Theta_{\mathcal{E}}}) $ and assign it to the noise-free data variable $\data$, and solve the inverse problem for the input $x$ from the perturbed code $\noisydata$
\begin{align*}
    \mathcal{E}(x,\mathbf{\Theta}_{\mathcal{E}}) = \noisydata \,\, .
\end{align*}

To be more precise, we first train a two-layer fully connected autoencoder $y = W_2(\sigma(W_1 x+b_1))+b_2$ using the vanilla stochastic gradient method (SGM) by minimising the mean squared error (MSE) on the MNIST training dataset. We set the code dimension to 100 and use ReLU as the activation function. Hence $\mathbf{\Theta}_{\mathcal{E}} = (W_1,b_1)$ where $W_1 \in \R^{784 \times 100}$ and $b_1 \in \R^{100 \times 1}$. 

All MNIST images are centred as a means of pre-processing. Algorithmically, we follow \eqref{eq:perceptron-inversion-tv-pdhg} to computationally solve \eqref{eq:generalised-perceptron-inversion}. The stepsize-parameters are chosen at $\tau_x = 1.99/\|W_1\|^2_2$ and $\tau_z = 1/(8\alpha)$. We choose the regularisation parameter $\alpha$ in the range $[10^{-4},10^{-2}]$ and set to $5 \times 10^{-3}$ for all sample images from the training set, and set to $\alpha = 5 \times 10^{-2}$ for all sample images from the validation set. These choices work well with regards to the visual quality of the inverted images. 

In Figure \ref{fig:perceptron-inversion-visualisation-train-images} and Figure \ref{fig:perceptron-inversion-visualisation-val-images}, we show visualisations of five sample images from the training set, and from the validation set respectively. In comparison, we have also visualised the decoder output. As can be seen, using the code that contains the same compressed information, the inverted images show more clearly defined edges and better visual quality than the decoded outputs. This is to be expected as we compare a nonlinear regularised inversion method with a linear decoder. 

\subsection{Multi-layer perceptrons}\label{sec:mlp-numerical-results}

\begin{figure}[!t]
\centering
\begin{minipage}{0.45\textwidth}
        \centering 
        \includegraphics[width = \textwidth]{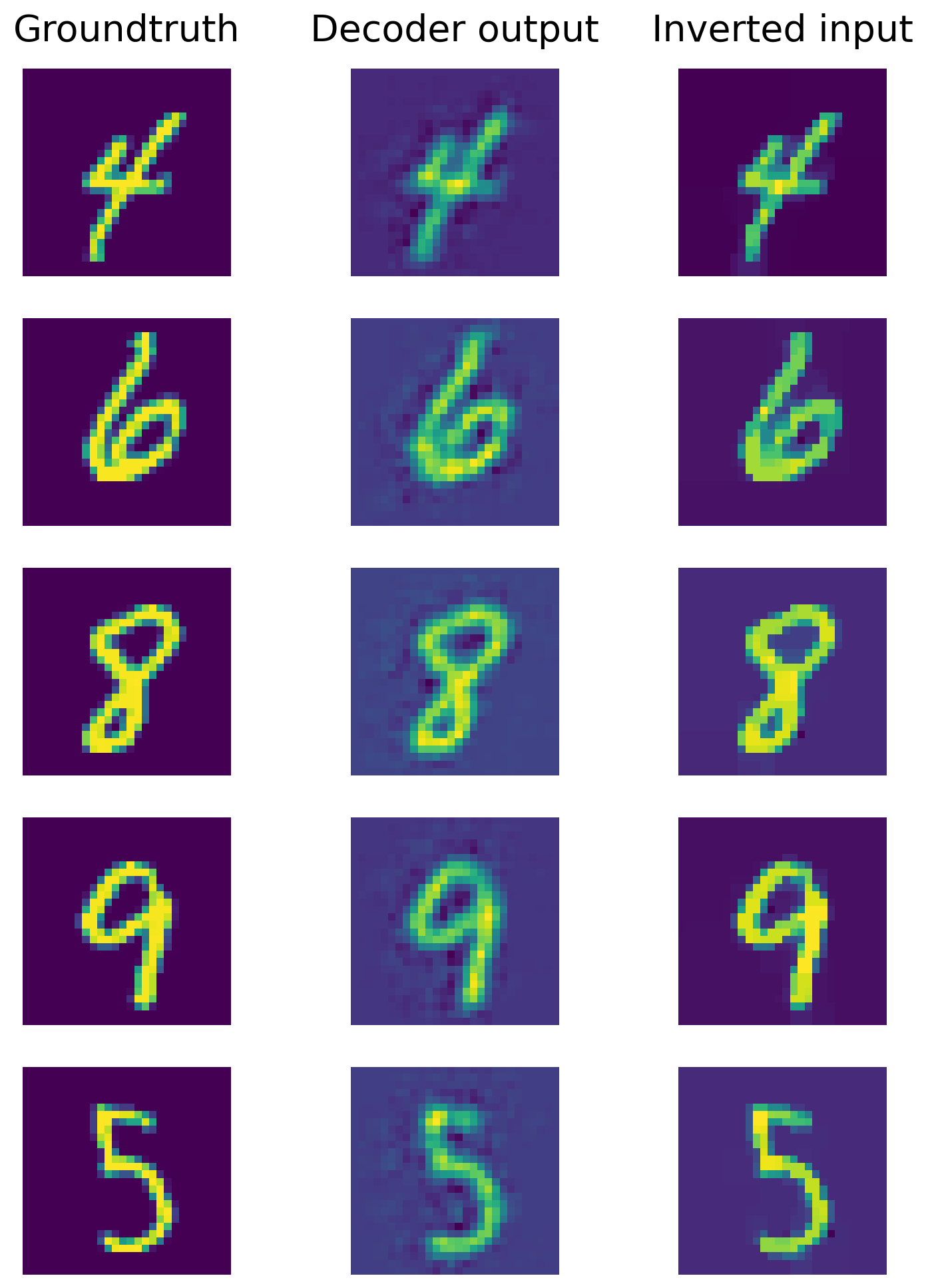}
        \caption{Inverted input images of the CNN from the MNIST training dataset, together with well-trained autoencoder output images and groudtruth input images}
        \label{fig:cnn-inversion-visualisation-train-images}
\end{minipage}
\begin{minipage}{0.45\textwidth}
        \centering
        \includegraphics[width = \textwidth]{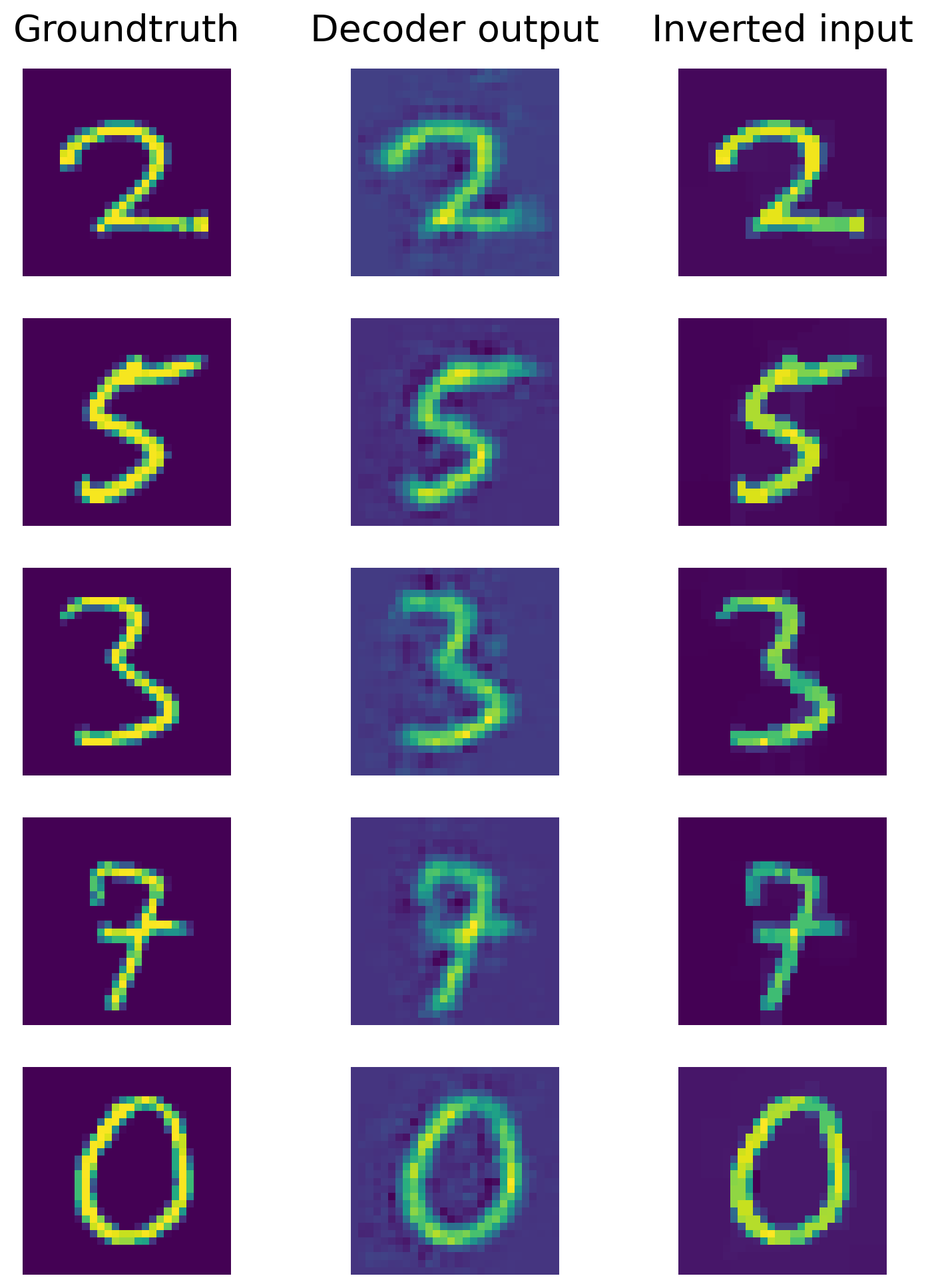}
        \caption{Inverted input images of the CNN from the MNIST validation dataset, together with well-trained autoencoder output images and groudtruth input images}
        \label{fig:cnn-inversion-visualisation-val-images}
\end{minipage}
\end{figure}


In this section, we present numerical results for inverting multi-layer perceptrons. In particular, we consider feedforward neural networks with convolutional layers (CNN), where in the network architecture two-dimensional convolution operations are used to represent the linear operations in the affine-linear functions $f(x,\Theta)$. Similar to the experimental design described in Section \ref{sec:perceptron-numerical-results}, we consider a multi-layer neural network inversion problem where we infer input image $x$ from a noise perturbed code $\noisydata$.

More specifically, we first train a six-layer convolutional autoencoder on the MNIST training dataset via stochastic gradient method to minimise the MSE. The encoder $\mathcal{E}(x,\mathbf{\Theta}_{\mathcal{E}})$ consists of two convolutional layers, both with $4 \times 4$ convolutions with stride $2$, each followed by the application of a ReLU activation function. As image spatial dimension reduce by half, we double the number of feature channels from 8 to 16. We use a fully-connected layer with weights $W_3 \in \R^{300 \times 784}$ and bias $b_3 \in \R^{300 \times 1}$ to generate the code. The decoder network first 
expands the code with an affine-linear transformation with weights $W_4 \in \R^{784 \times 300}$ and bias $b_4 \in \R^{784 \times 1}$. This is followed by two layers of transpose convolutionals with kernel size $4 \times 4$, where each is followed by a ReLU activation function. The number of feature channels halves each time as we double the spatial dimension. 

Following the implementation details outlined in Section \ref{sec:implementation-mlp}, we iteratively compute the update steps \eqref{eq:mlp-x0-pdhg} and \eqref{eq:x_l-proximal-gradient-specific-form}. For the PDHG method, we choose the stepsize-parameters as $\tau_x = 1.99/\|W_1\|^2_2$ and $\tau_z = 1/(8\alpha)$. The initial values $x_0$ and $z$ are both zero. The update steps stop either after reaching the maximum iterations of $1500$ or when the improvements on $x_0$ and $z$ are less than $10^{-5}$ in norm. For the coordinate descent algorithm, the stepsize-parameters are set to $\tau_{x_{l}} = 1.99/\|W_{l+1}\|^2_2$ for each layer. 

In Figure \ref{fig:cnn-inversion-visualisation-train-images} and Figure \ref{fig:cnn-inversion-visualisation-val-images}, we visualise the inverted images, the decoder output images, along with the groundtruth images, from the training dataset and validation dataset respectively. For each image, $\alpha$ is chosen in the range $[10^{-4},10^{-2}]$ and set at $9\times 10^{-3}$ for both training sample images and validation sample images for best visual inversion quality. 

\begin{figure}
    \centering
    \includegraphics[scale=0.5]{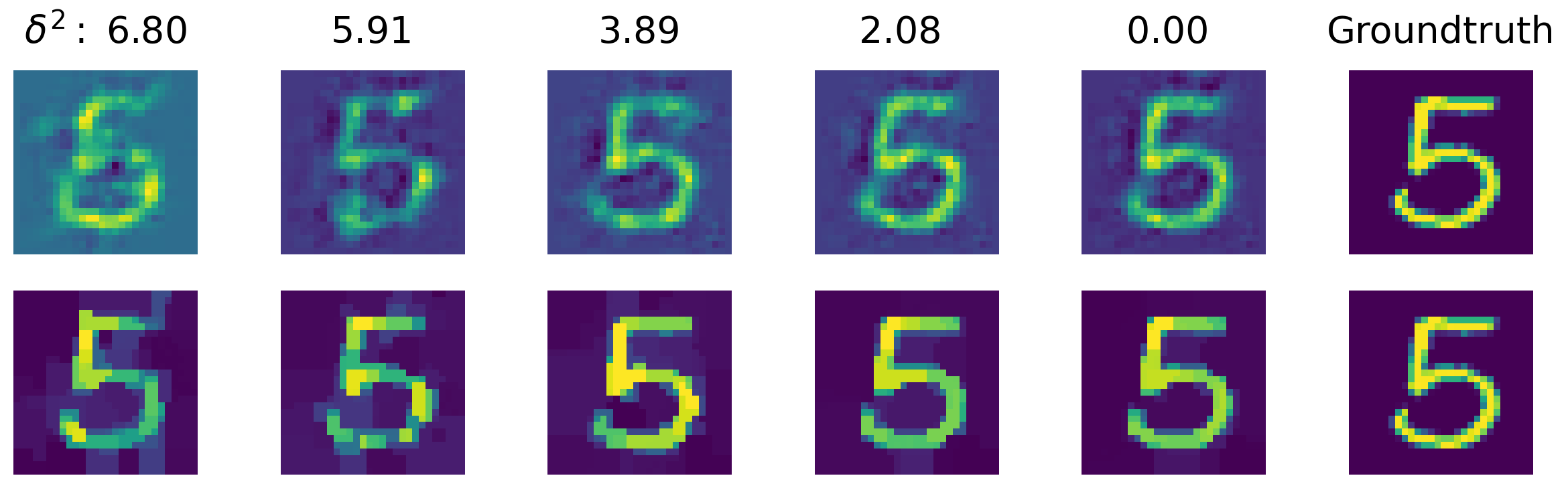}
    \caption{Visualisation of the comparison between inverted image and decoded image against various levels of noise. \textbf{Top: } Decoded output image from the trained convolutional autoencoder. \textbf{Bottom: } Inverted input image from the CNN with total variation regularisation.}
    \label{fig:noisy_comparison}
\end{figure}

In Figure \ref{fig:noisy_comparison} we further compare how total variation regularisation and decoder respond to different levels of data noise. The noisy data is produced by adding Gaussian noise to perturb the code of each image. We start with zero mean Gaussian noise with standard deviation 0.33 and gradually reduce the noise level, this translates to decreasing $\delta^2$ from 6.80 down to 0.00. 

Please note that for each noise level the regularisation factor $\alpha$ is manually selected in the range $[10^{-4}, 10^{-2}]$ for the best PSNR value. As we can see, for the noise level with standard deviation $0.33$ where $\delta^2$ is at 6.80, the decoder is only capable of producing a blurry distorted output, while the inverted image shows the structure of the digit more clearly. When we decrease the noise level down to $0.00$, the inverted image becomes more clean-cut while the decoded image is still less sharply defined.

\begin{figure}
    \centering
    \includegraphics[scale=0.4]{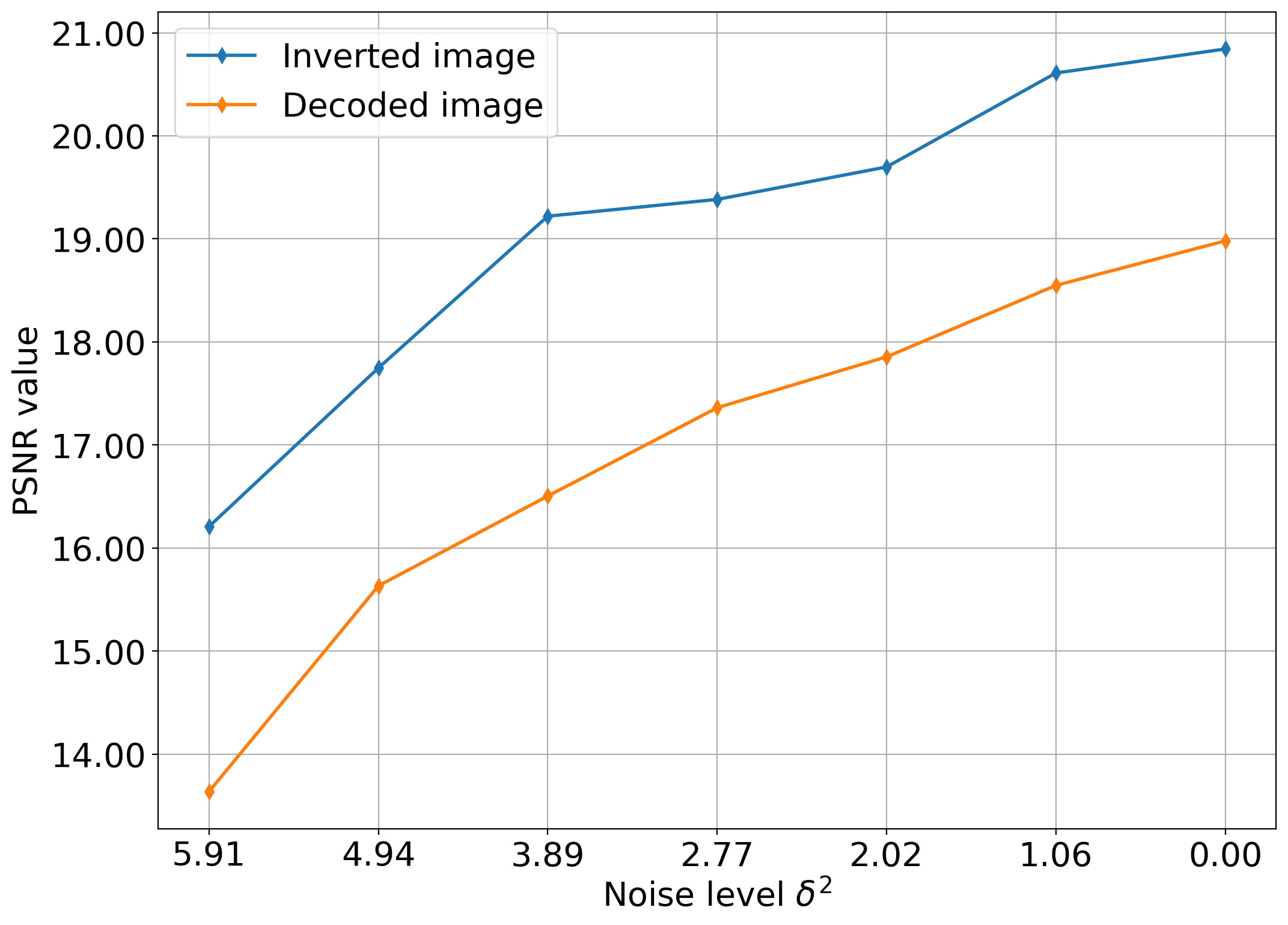}
    \caption{Comparison of PSNR values of total variation-based reconstruction and decoder output per noise level. Each curve reports the change of PSNR value over gradually decreasing levels of Gaussian noise, with $\delta^2$ ranging from 0.00 to 6.80.} 
    \label{fig:psnr_plot}
\end{figure}

Figure \ref{fig:psnr_plot} plots the PSNR value of the decoded and inverted image against decreasing noise level. We want to emphasise that it would be more rigorous to compute and compare $\symbreg(\reconstruction[\alpha(\noisebound)], \solution)$ as suggested in the error estimate bound in \eqref{eq:perceptron-error-estimate}, but empirically the PSNR value does also support the notion of a convergent regularisation. 

\section{Conclusions \& Outlook}\label{sec:outlook}
We have introduced a novel variational regularisation framework based on a lifted Bregman formulation for the stable inversion of feed-forward neural networks (also known as multi-layer perceptrons). We have proven that the proposed framework is a convergent regularisation for the single-layer perceptron case under the mild assumption that the inverse problem solution has to be in the range of the regularisation operator. We have derived a general error estimate as well as a specific error estimate for the case that the activation function is the ReLU activation function. We have also addressed the extension of the theory to the multi-layer perceptron case, which can be carried out sequentially, albeit under unrealistic assumptions. We have discussed implementation strategies to solve the proposed scheme computationally, and presented numerical results for the regularised inversion of the image of a circle and piecewise constant images of hand-written digits from single- and multi-layer perceptron outputs with total variation regularisation.

Despite all the positive achievements presented in this work, the proposed framework also has some limitations. The framework is currently restricted to feed-forward architectures with affine-linear transformations and proximal activation functions. While it is straight-forward to extend the framework to other architectures such as ResNets \citep{he2016deep} or U-Nets \citep{ronneberger2015u}, it is not straight-forward to include nonlinear operations that cannot be expressed as proximal maps of convex functions, such as max-pooling. However, for many examples there exist remedies, such as using average pooling instead of max-pooling in the previous example. 

An open question is how a convergence theory without restrictive, unrealistic assumptions can be established for the multi-layer case. One issue is the non-convexity of the proposed formulation. A remedy could be the use of different architectures that lead to lifted Bregman formulations that are jointly convex in all auxiliary variables. 

And last but not least, one would also like to consider other forms of regularisation, such as iterative regularisation, data-driven regularisations \citep{kabri2022convergent}, or even combinations of both \citep{data_driven_landweber}. However, a convergence analysis for such approaches is currently an open problem.


\section*{Conflict of Interest Statement}
The authors declare that the research was conducted in the absence of any commercial or financial relationships that could be construed as a potential conflict of interest.

\section*{Author Contributions}


XW has programmed and contributed all numerical results as well as Section \ref{sec:implementation} and Section \ref{sec:results}. MB has contributed the introduction (Section \ref{sec:introduction}) as well as the theoretical results (Section \ref{sec:convergence}). Both authors have contributed equally to Section \ref{sec:feed-forward} and Section \ref{sec:outlook}.


\section*{Acknowledgments}
The authors acknowledge support from the Cantab Capital Institute for the Mathematics of Information, the Cambridge Centre for Analysis (CCA) and the Alan Turing Institute (ATI).


\section*{Data Availability Statement}
The programming code for this study can be found in the University
of Cambridge data repository at 
10.17863/CAM.94404.

\bibliographystyle{Frontiers-Harvard} 
\bibliography{main}

\end{document}